\input amssym.def
\input epsf

\let \blskip = \baselineskip
\parskip=1.2ex plus .2ex minus .1ex

\tabskip 20pt
\tolerance = 1000
\pretolerance = 50
\newcount\itemnum
\itemnum = 0
\overfullrule = 0pt

\def\title#1{\bigskip\centerline{\bigbigbf#1}}
\def\author#1{\bigskip\centerline{\bf #1}\smallskip}
\def\address#1{\centerline{\it#1}}
\def\abstract#1{\vskip1truecm{\narrower\noindent{\bf Abstract.} #1\bigskip}}

\def\sp{\bigskip}
\def\nosp{\vskip -\the\blskip plus 1pt minus 1pt}

\def\br{\hfil\break} 
\def\ti{\br \hglue \the \parindent}

\def\skipit#1{}
\def\mdag{\raise 3pt\hbox{\dag}}

\def\XP{\par\noindent\hang}
\def\LP{\par\noindent}
\def\BP[#1]{\par\item{[#1]}}
\def\SH#1{\sp\vskip\parskip\leftline{\bigbf #1}\nobreak}

\def\TH#1{\sp\XP{\bf THEOREM\ \shead#1}}
\def\LM#1{\sp\XP{\bf LEMMA\ \shead#1}}

\def\PR#1{\sp\XP{\bf PROPOSITION\ \shead#1}}
\def\CO#1{\sp\XP{\bf COROLLARY\ \shead#1}}

\def\EX#1{\sp\LP{\bf Example\ \shead#1}}
\def\PF{\LP{\bf Proof:\ }}

\def\NX{\advance\itemnum by 1 \sp\LP {\bf \shead \the\itemnum.\ }}
\def\qed{\null\nobreak\hfill\hbox{${\vrule width 5pt height 6pt}$}\par\sp}

\def\cart{\>\hbox{${\vcenter{\vbox{
    \hrule height 0.4pt\hbox{\vrule width 0.4pt height 4.5pt
    \kern4pt\vrule width 0.4pt}\hrule height 0.4pt}}}$}\>}
\def\bxmu{\>\hbox{${\vcenter{\vbox {
    \hrule height 0.4pt\hbox{\vrule width 0.4pt height 4pt
    \hskip -1.3pt\lower 1.8pt\hbox{$\times$}\negthinspace\vrule width 0.4pt}
    \hrule height 0.4pt}}}$}\>}

\def\lin#1{\hbox to #1true in{\hrulefill}}



\def\JGT{{\it J.\ Graph Theory}}

\def\DM{{\it Discrete Math.{}}}

\def\ADM{{\it Annals Discr.\ Math.{}}}

\def\SIAD{{\it SIAM\ J.\ Algeb.\ Disc.\ Meth.{}}}
\def\SIAP{{\it SIAM\ J.\ Appl.\ Math.{}}}
\def\SIDM{{\it SIAM\ J.\ Discr.\ Math.{}}}

\def\al{\alpha}			\def\dlt{\delta}
\def\eps{\epsilon}    
	 	\def\DLT{\Delta}

 \def\ZZ{{\Bbb Z}}


\def\esub{\subseteq}



\def\({\left(}	\def\){\right)}


\def\CH#1#2{{{#1}\choose{#2}}} 
\def\FR#1#2{{#1 \over #2}}     

\def\FL#1{\left\lfloor{#1}\right\rfloor} \def\FFR#1#2{\FL{{#1\over#2}}}
\def\CL#1{\left\lceil{#1}\right\rceil}   \def\CFR#1#2{\CL{{#1\over#2}}}

\def\VEC#1#2#3{#1_{#2},\ldots,#1_{#3}}
\def\VECOP#1#2#3#4{#1_{#2}#4\cdots #4 #1_{#3}}

\def\MAP#1#2#3{#1\colon\,#2\to#3}
\def\SET#1:#2{\{#1\colon\;#2\}}

		
\def\C#1{\left | #1 \right |}    


    \def\diam{ {\rm diam\;}}


\magnification=\magstep1
\vsize=9.0 true in
\hsize=6.5 true in
\headline={\hfil\ifnum\pageno=1\else\folio\fi\hfil}
\footline={\hfil\ifnum\pageno=1\folio\else\fi\hfil}

\parindent=20pt
\baselineskip=12pt
\parskip=.2ex  

\def\shead{ }

\font\bigbf = cmb10 scaled \magstep1

\font\bigbigbf = cmb10 scaled \magstep2




\title{EDGE-BANDWIDTH OF GRAPHS}
\author{Tao Jiang}
\address{University of Illinois, Urbana, IL 61801-2975}
\author{Dhruv Mubayi}
\address{Georgia Institute of Technology, Atlanta, GA 30332-0160}
\author{Aditya Shastri}
\address{Banasthali University, Rajasthan 304 022, India}
\author{Douglas B. West\mdag}
\address{University of Illinois, Urbana, IL 61801-2975, {west@math.uiuc.edu}}
\vfootnote{}{\br
   \mdag Research supported in part by NSA/MSP Grant MDA904-93-H-3040.\br
   Running head: EDGE-BANDWIDTH \br
}
\abstract{
The {\it edge-bandwidth} of a graph is the minimum, over all labelings
of the edges with distinct integers, of the maximum difference between
labels of two incident edges.  We prove that edge-bandwidth is at least
as large as bandwidth for every graph, with equality for certain caterpillars.
We obtain sharp or nearly-sharp bounds on the change in edge-bandwidth under
addition, subdivision, or contraction of edges.  We compute edge-bandwidth for
$K_n$, $K_{n,n}$, caterpillars, and some theta graphs.}

\SH
{1. INTRODUCTION}

A classical optimization problem is to label the vertices of a graph with
distinct integers so that the maximum difference between labels on adjacent
vertices is minimized.  For a graph $G$, the optimal bound on the differences is
the {\it bandwidth} $B(G)$.  The name arises from computations with sparse
symmetric matrices, where operations run faster when the matrix is permuted so
that all entries lie near the diagonal.  The bandwidth of a matrix $M$ is the
bandwidth of the corresponding graph whose adjacency matrix has a 1 in those
positions where $M$ is nonzero.  Early results on bandwidth are surveyed in [2]
and [3].

In this paper, we introduce an analogous parameter for edge-labelings.
An {\it edge-numbering} (or {\it edge-labeling}) of a graph $G$ is a function
$f$ that assigns distinct integers to the edges of $G$.  We let $B'(f)$ denote
the maximum of the difference between labels assigned to adjacent (incident)
edges.  The {\it edge-bandwidth} $B'(G)$ is the minimum of $B'(f)$ over all
edge-labelings.  The term ``edge-numbering'' is used because we may assume that
$f$ is a bijection from $E(G)$ to the first $|E(G)|$ natural numbers.

We use the notation $B'(G)$ for the edge-bandwidth of $G$ because it is
immediate that the edge-bandwidth of a graph equals the bandwidth of its
line graph.  Thus well-known elementary bounds on bandwidth can be applied
to line graphs to obtain bounds on edge-bandwidth.  We mention several
such bounds.  We compute edge-bandwidth on a special class where all these
bounds are arbitrarily bad.

The relationship between edge-bandwidth and bandwidth is particularly
interesting.  Always $B(G)\le B'(G)$, with equality for caterpillars of
diameter more than $k$ in which every vertex has degree 1 or $k+1$.  Among
forests, $B'(G)\le 2B(G)$, which is almost sharp for stars.  More generally, if
$G$ is a union of $t$ forests, then $B'(G)\le 2tB(G)+t-1$.

Chv\'atalov\'a and Opatrny [5] studied the effect on bandwidth of edge addition,
contraction, and subdivision (see [22] for further results on edge addition).
We study these for edge-bandwidth.  Adding or contracting an edge at most
doubles the edge-bandwidth.  Subdividing an edge decreases the edge-bandwidth by
at most a factor of $1/3$.  All these bounds are sharp within additive
constants.  Surprisingly, subdivision can also increase edge-bandwidth, but
at most by 1, and contraction can decrease it by 1.

Because the edge-bandwidth problem is a restriction of the bandwidth problem,
it may be easier computationally.  Computation of bandwidth is NP-complete [17],
remaining so for trees with maximum degree 4 [8] and for several classes of
caterpillar-like graphs.  Such graphs generally are not line graphs (they
contain claws).  It remains open whether computing edge-bandwidth
(computing bandwidth of line graphs) is NP-hard.

Due to the computational difficulty, bandwidth has been studied on various
special classes.  Bandwidth has been determined for caterpillars and for various
generalizations of caterpillars ([1,11,14,21]), for complete $k$-ary trees [19],
for rectangular and triangular grids [4,10] (higher dimensions [9,15]), for
unions of pairwise internally-disjoint paths with common endpoints (called
``theta graphs'' [6,13,18]), etc.  Polynomial-time algorithms exist for
computing bandwidth for graphs in these classes and for interval graphs [12,20].
We begin analogous investigations for edge-bandwidth by computing the
edge-bandwidth for cliques, for equipartite complete bipartite graphs, and for
some theta graphs.

\sp
\SH
{2. RELATION TO OTHER PARAMETERS}

We begin by listing elementary lower bounds on edge-bandwidth that follow from
standard arguments about bandwidth when applied to line graphs.

\PR{1.} Edge-bandwidth satisfies the following.
\br
a) $B'(H)\le B'(G)$ when $H$ is a subgraph of $G$.
\br
b) $B'(G) = \max \{B'(G_i)\}$, where $\{G_i\}$ are the components of $G$.
\br
c) $B'(G)\ge\DLT(G)-1$.
\PF
(a) A labeling of $G$ contains a labeling of $H$.  (b) Concatenating labelings
of the components achieves the lower bound established by (a).  (c) The edges
incident to a single vertex induce a clique in the line graph.  The lowest and
highest among these labels are at least $\DLT(G)-1$ apart.  \qed

\vskip -1pc
\PR{2.}
$B'(G)\ge\max_{H\esub G}\CFR{e(H)-1}{\diam(L(H))}$.
\PF
This is the statement of Chung's ``density bound'' [3] for line graphs.
Every labeling of a graph contains a labeling of every subgraph.
In a subgraph $H$, the lowest and highest labels are at least $e(H)-1$
apart, and the edges receiving these labels are connected by a path
of length at most $\diam(L(H))$, so by the pigeonhole principle some
consecutive pair of edges along the path have labels differing by at least
$(e(H)-1)/\diam(L(H))$.  \qed

Subgraphs of diameter 2 include stars, and a star in a line graph is
generated from an edge of $G$ with its incident edges at both endpoints.
The size of such a subgraph is at most $d(u)+d(v)-1$, yielding the bound
$B'(G)\ge [d(u)+d(v)]/2-1$ for $uv\in E(G)$.  This is at most $\DLT(G)-1$, the
lower bound from Proposition 1.  Nevertheless, because of the way in which
stars in line graphs arise, they can yield a better lower bound for
regular or nearly-regular graphs.  We develop this next.

\PR{3.}
For $F\esub E(G)$, let $\partial(F)$ denote the set of edges not in $F$ that
are incident to at least one edge in $F$.  The edge-bandwidth satisfies
$B'(G)\ge\max_k\min_{|F|=k}|\partial(F)|$.
\PF
This is the statement of Harper's ``boundary bound'' [9] for line graphs.
Some set $F$ of $k$ edges must be the set given the $k$ smallest labels.
If $m$ edges outside this set have incidences with this set, then the largest
label on the edges of $\partial F$ is at least $k+m$, and the difference between
the labels on this and its incident edge in $F$ is at least $m$.  \qed

\vskip -.5pc
\CO{4.}
$B'(G)\ge\min_{uv\in E(G)} d(u)+d(v)-2$.
\PF
We apply Proposition 3 with $k=1$.  Each edge $uv$ is incident to $d(u)+d(v)-2$
other edges.  Some edge must have the least label, and this establishes
the lower bound.  \qed

Although these bounds are often useful, they can be arbitrarily bad.
The {\it theta graph} $\Theta(\VEC l1m)$ is the graph that is the union of $m$
pairwise internally-disjoint paths with common endpoints and lengths $\VEC l1m$.
The name ``theta graph'' comes from the case $m=3$.  The bandwidth is known for
all theta graphs, but settling this was a difficult process finished in
[18].  When the path lengths are equal, the edge-bandwidth and bandwidth both
equal $m$, using the density lower bound and a simple construction.  The
edge-bandwidth can be much higher when the lengths are unequal.  Our example
showing this will later demonstrate sharpness of some bounds.

Our original proof of the lower bound was lengthy.  The simple argument
presented here originated with Dennis Eichhorn and Kevin O'Bryant.  It will be
generalized in [7] to compute edge-bandwidth for a large class of theta graphs.

\EX{A.}
Consider $G=\Theta(\VEC l1m)$ with $l_m=1$ and $\VECOP l1{m-1}= = 3$.  Let
$a_i,b_i,c_i$ denote the edges of the $i$th path of length 3, and let $e$ be the
edge incident to all $a_i$'s at one end and to all $c_i$'s at the other end.
Since $\DLT(G)=m$, Proposition 1c yields $B'(G)\ge m-1$.  Proposition 2 also
yields $B'(G)\ge m-1$.  For $1\le k\le 2m-2$, the first $k$ edges in the list
$\VEC a1{m-1},\VEC b1{m-1}$ are together incident to exactly $m$ other edges,
and larger sets are incident to at most $m-1$ other edges.  Thus the best lower
bound from Proposition 3 is at most $m$.

Nevertheless, $B'(G)=\CL{(3m-3)/2}$.  
For the upper bound, we assign the $3m-2$ labels in order to $a$'s, $b$'s, and
$c$'s, inserting $e$ before $b_{\CL{m/2}}$.
The difference between labels of incidence edges is always at most $m$ except
for incidences involving $e$, which are at most $\CL{(3m-3)/2}$ since $e$
has the middle label.
$$\VEC a1{m-1},\VEC b1{\CL{m/2}-1},e,\VEC b{\CL{m/2}}{m-1},\VEC c1{m-1}.$$

To prove the lower bound, consider a numbering $\MAP f{E(G)}{\ZZ}$, and let
$k=B'(f)$.  Let $\al=\max\{f(e),\max_i\{f(a_i)\}\}$ and
$\al'=\min\{f(e),\min_i\{f(c_i)\}\}$.  Comparing the edges with labels
$\al,f(e),\al'$ yields $\al-k\le f(e)\le\al'+k$.  Let $I$ be the interval
$[\al-k,\al'+k]$.  By construction, $I$ contains the labels of all $a$'s, all
$c$'s, and $e$.  If $f(a_i)<\al'$ and $f(c_i)>\al$, then also $f(b_i)\in I$.
By the choice of $\al,\al'$, avoiding this requires $\al'<f(a_i)\le\al$ or
$\al'\le f(c_i)<\al$.  Since each label is assigned only once and the label
$f(e)$ cannot play this role, only $\al-\al'$ of the $b$'s can have labels
outside $I$.  Counting the labels we have forced into $I$ yields
$\C I\ge (2m-1)+(m-1-\al+\al')$.  On the other hand, $\C I = 2k+\al'-\al+1$.
Thus $k\ge (3m-3)/2$, as desired.
\qed

\sp
\SH
{3. EDGE-BANDWIDTH VS. BANDWIDTH}
In this section we prove various best-possible inequalities involving
bandwidth and edge-bandwidth.  The proof that $B(G)\le B'(G)$ requires
several steps.  All steps are constructive.
When $f$ or $g$ is a labeling of the edges or vertices of $G$, we say that
$f(e)$ of $g(v)$ is the $f$-{\it label} or $g$-{\it label} of the edge $e$
or vertex $v$. An $f$-{\it label} on an edge incident to $u$ is an
{\it incident $f$-label} of $u$.

\LM{5.}
If a finite graph $G$ has minimum degree at least two, then $B(G)\le B'(G)$.
\PF
From an optimal edge-numbering $f$ (such that $B'(f)=B'(G)=m$), we define a
labeling $g$ of the vertices.  The labels used by $g$ need not be consecutive,
but we show that $|g(u)-g(v)|\le m$ when $u$ and $v$ are adjacent.

We produce $g$ in phases.  At the beginning of each phase, we choose
an arbitrary unlabeled vertex $u$ and call it the {\it active vertex}.
At each step in a phase, we select the unused edge $e$ of smallest $f$-label
among those incident to the active vertex.  We let $f(e)$ be the $g$-label
of the active vertex, mark $e$ {\it used}, and designate the other endpoint
of $e$ as the active vertex.  If the new active vertex already has a label,
we end the phase.  Otherwise, we continue the phase.

When we examine a new active vertex, it has an edge with least incident 
label, because every vertex has degree at least 2 and we have not previously
reached this vertex.  Each phase eventually ends, because the vertex set is
finite and we cannot continue reaching new vertices.  The procedure assigns a
label $g(u)$ for each $u\in V(G)$, since we continue to a new phase as long as
an unlabeled vertex remains.

It remains to verify that $|g(u)-g(v)|\le m$ when $uv\in E(G)$.
Suppose that $g(u)=a=f(e)$ and $g(v)=b=f(e')$.  Since each vertex is assigned
the $f$-label of an incident edge, we have $e,e'$ incident to $u,v$,
respectively.  If the edge $uv$ is one of $e,e'$, then $e$ and $e'$ are
incident, which implies that $|g(u)-g(v)|=|f(e)-f(e')|\le m$.

Otherwise, we have $f(uv)=c$ for some other value $c$.  We may assume that
$a<b$ by symmetry.  If $a<c$ and $b<c$, then
$|g(u)-g(v)|=b-a<c-a=f(uv)-f(e)\le m$.
Thus we may assume that $b>c$.  In particular, $g(v)$ is not the least $f$-label
incident to $v$.

The algorithm assigns $v$ a label when $v$ first becomes active, using the least
$f$-label among {\it unused} incident edges.  When $v$ first becomes
active, only the edge of arrival is a used incident edge.  Thus $g(v)$ is the
least incident $f$-label except when $v$ is first reached via the least-labeled
incident edge.  In this case, $g(v)$ is the second smallest incident $f$-label.
Thus $c$ is the least $f$-label incident to $v$ and $v$ becomes active by
arrival from $u$.  This requires $g(u)=c$, which contradicts $g(u)=a$ and
eliminates the bad case.  \qed

\vskip -.5pc
\LM{6.}
If $G$ is a tree, then $B(G)\le B'(G)$.
\PF
Again we use an optimal edge-numbering $f$ to define a vertex-labeling
$g$ whose adjacent vertices differ by at most $B'(f)$.  We may assume that
the least $f$-label is 1, occurring on the edge $e=uv$.  Assign (temporarily)
$g(u)=g(v)=f(e)$.  View the edge $e$ as the root of $G$.  For each vertex
$x\notin\{u,v\}$, let $g(x)$ be the $f$-label of the edge incident to $x$
along the path from $x$ to the root.

If $xy\in E(G)$ and $xy\ne uv$, then we may assume that $y$ is on the path from
$x$ to the root.  We have assigned $g(x)=f(xy)$, and $g(y)$ is the $f$-label
of an edge incident to $y$, so $|g(x)-g(y)|\le B'(f)$.

Our labeling $g$ fails to be the desired labeling only because we used 1
on both $u$ and $v$.  Observe that the largest $f$-label incident to $uv$
occurs on an edge incident to $u$ or on an edge incident to $v$ but not
both; we may assume the latter.  Now we change $g(u)$ to 0.  Because the
differences between $f(uv)$ and $f$-labels on edges incident to $u$
were {\it less} than $B'(f)$, this produces the desired labeling $g$.  \qed

\vskip -.5pc
\TH{7.}
For every graph $G$, $B(G)\le B'(G)$.
\PF
By Proposition 1b, it suffices to consider connected graphs.
Let $f$ be an optimal edge-numbering of $G$; we produce a vertex labeling
$g$.  Lemma 6 applies when $G$ is a tree.   Otherwise, $G$ contains a
cycle, and iteratively deleting vertices of degree 1 produces a
subgraph $G'$ in which every vertex has degree at least 2.
The algorithm of Lemma 5, applied to the restriction of $f$ to $G'$,
produces a vertex labeling $g$ of $G'$ in which (1) adjacent vertices have
labels differing by at most $B'(f)$, and (2) the label on each
vertex is the $f$-label of some edge incident to it in $G'$.

To obtain a vertex labeling of $G$, reverse the deletion procedure.  This
iteratively adds a vertex $x$ adjacent to a vertex $y$ that already has a
$g$-label.  Assign to $x$ the $f$-label of the edge $xy$ in the full
edge-numbering $f$ of $G$.  Now $g(x)$ and $g(y)$ are the $f$-labels of two
edges incident to $y$ in $G$, and thus $|g(x)-g(y)|\le B'(f)$.  The claims
(1) and (2) are preserved, and we continue this process until we replace all
vertices that were deleted from $G$.
\qed

A {\it caterpillar} is a tree in which the subtree obtained by deleting all
leaves is a path.  One of the characterizations of caterpillars is the
existence of a linear ordering of the edges such that each prefix and
each suffix forms a subtree.  We show that such an ordering is optimal
for edge-bandwidth and use this to show that Theorem 7 is nearly sharp.

\PR{8.}
If $G$ is a caterpillar, then $B'(G)=\DLT(G)-1$.
Let $G$ be the caterpillar of diameter $d$ in which every vertex
has degree $k+1$ or 1.  If $d\ge k$, then $B(G)=B'(G)=k$.
\PF
Let $G$ be a caterpillar.  Let $\VEC v1{d-1}$ be the non-leaf vertices
of the dominating path.  The diameter of $G$ is $d$.  Number the edges
by assigning labels in the following order: first the pendant edges
incident to $v_1$, then $v_1v_2$, then the pendant edges incident to $v_2$,
then $v_2v_3$, etc.  Since edges are incident only at $\VEC v1{d-1}$,
this ordering places all pairs of incident edges within $\DLT(G)-1$
positions of each other.  Since $B'(G)\ge \DLT(G)-1$ for all $G$,
equality holds.

For a caterpillar $G$ with order $n$ and diameter $d$,
Chung's density bound yields $B(G)\ge (n-1)/d$.
Let $G$ be the caterpillar of diameter $d$ in which every vertex
has degree $k+1$ or 1.  We have $d-1$ vertices of degree $k+1$,
so $n=(d-1)k+2$ and $B(G)> k-k/d$.  When $d\ge k$, we have $B(G)\ge k$.

On the other hand, we have observed that $B'(G)\le \DLT(G)-1 =k$
for caterpillars.  By Theorem 7, equality holds throughout for these
special caterpillars.  \qed

Theorem 7 places a lower bound on $B'(G)$ in terms of $B(G)$.
We next establish an upper bound.  The {\it arboricity} 
is the minimum number of forests needed to partition the edges of $G$.

\TH{9.}
If $G$ has arboricity $t$, then $B'(G)\le 2tB(G)+t-1$.
When $t=1$, the inequality is almost sharp; there are caterpillars with
$B'(G)=2B(G)-1$.
\PF
Given an optimal number $g$ of $V(G)$, we construct a labeling $f$ of $E(G)$.
Let $\VEC G1t$ be a decomposition of $G$ into the minimum number of forests.
In each component of each $G_i$, select a root.  Each edge of $G_i$ is the first
edge on the path from one of its endpoints to the root of its component in
$G_i$; for $e\in E(G_i)$, let $v(e)$ denote this endpoint.
Define $f(e)=tg(v(e))+i$.

Each vertex of each forest heads toward the root of its component in that forest
along exactly one edge, so the $f$-labels of the edges are distinct.
Each $f$-label arises from the $g$-label of one of its endpoints.
Thus the $f$-labels of two incident edges arise from the $g$-labels of
vertices separated by distance at most 2 in $G$.  Also the indices of the
forests containing these edges differ by at most $t-1$.
Thus when $e,e'$ are incident we have $|f(e)-f(e')|\le t2B(g)+t-1$.

The bandwidth of a caterpillar is the maximum density (\#edges/diameter) over
subtrees [14].  This equals $\CL{\DLT(G)/2}$ whenever the vertex degrees all
lie in $\{\DLT(G),2,1\}$ and the vertices of degree $\DLT(G)$ are pairwise
nonadjacent.  When $\DLT(G)$ is even, Proposition 8 yields $B'(G)=2B(G)-1$.
(Without [14], this still holds explicitly for stars.)
\qed

\SH
{4. EFFECT OF EDGE OPERATIONS}

In this section, we obtain bounds on the effect of local edge operations on
the edge-bandwidth.  The variations can be linear in the value of the
edge-bandwidth, and our bounds are optimal except for additive constants.
We study addition, subdivision, and contraction of edges.

\TH{10.}
If $H$ is obtained from $G$ by adding an edge, then
$B'(G)\le B'(H)\le 2B'(G)$.  Furthermore, for odd $k$ there are examples of
$H=G+e$ such that $B'(G)=k$ and $B'(H)\ge 2k-1$.
\PF
The first inequality holds because $G$ is a subgraph of $H$.
For the second, let $g$ be an optimal edge-numbering of $G$;
we produce an edge-numbering $f$ of $H$ such that $B'(f)\le2B'(g)$.

If $e$ is not incident to an edge of $G$, form $f$ from $g$ by giving $e$ a
new label higher than the others.  If only one endpoint of $e$ is incident to
an edge $e'$ of $G$, form $f$ by leaving the $g$-labels less than $g(e')$
unchanged, augmenting the remaining labels by 1, and letting $f(e)=g(e')+1$.
We have $B(f)\le B(g)+1$.

Thus we may assume that the new edge $e$ joins two vertices of $G$.
Our construction for this case modifies an argument in [22].
Let $e_i$ be the edge such that $g(e_i)=i$, for $1\le i\le B(g)$.
Let $p,q$ be the smallest and largest indices of edges of $G$ incident
to $e$, respectively, and let $r=\FL{(p+q)/2}$.

The idea in defining $f$ from $g$ is to ``fold'' the ordering at $r$,
renumbering out from there so that $e_p$ and $e_q$ receive consecutive
labels, and inserting $e$ just before this.  The renumbering of the old
edges is as follows
$$f(e_j)=\cases{2(j-r)&if $r<j<q$\cr
                2(j-r)+1&if $q\le j$\cr
                2(r-j)+1&if $p< j\le r$\cr
                2(r-j)+2&if $j\le p$\cr}$$
Finally, let $f(e)=\min\{f(e_p),f(e_q)\}-1=q-p$.  After the edges with
$g$-labels higher than $q$ or lower than $p$ are exhausted, the new numbering
leaves gaps.  For edges $e_i,e_j\in E(G)$,
we have $|f(e_i)-f(e_j)|\le 2|i-j|+1$, where the
possible added 1 stems from the insertion of $e$.  When $r$ is between $i$
and $j$, the actual stretch is smaller.

It remains to consider incidences involving $e$.  Suppose that $e'=e_j$
is incident to $e$.  Note that $1\le f(e')\le q-p+2 = f(e)+2$;
we may assume that $1\le f(e')< f(e)$.  If $e_p$ and $e_q$ are incident to the
same endpoint of $e$, then $1\le f(e)-f(e')\le q-p+1\le B(g)+1$.  If $e_p$ and
$e_q$ are incident to opposite endpoints of $e$, then $e'$ is incident
to $e_p$ or $e_q$.  In these two cases, we have $p\le j\le p+B(g)$ or 
$q-B(g)\le j\le q$.  Since $j$ differs from $p$ or $q$, respectively, by
at most $B(g)$, we obtain $1\le f(e)-f(e')\le 2B(g)$.

The bound is nearly sharp when $k$ is odd.  Let $G$ be the caterpillar of
diameter $k+1$ with vertices of degree $k+1$ and 1 (see Proposition 8).
We have $e(G)=k^2+1$ and $B'(G)=B(G)=k$.  The graph $H$ formed by
adding the edge $v_1v_k$ is a cycle of length $k$ plus pendant edges; each
vertex of the cycle has degree $k+1$ except for two adjacent vertices of degree
$k+2$.  The diameter of $L(H)$ is $\FL{k/2}+1=(k+1)/2$, and $H$ has $k^2+2$
edges.  By Proposition 2, we obtain
$B'(H)\ge \CFR{k^2+1}{(k+1)/2}=\CL{2k-2+\FR4{k+1}}=2k-1$.
\qed

{\it Subdividing} an edge $uv$ means replacing $uv$ by a path $u,w,v$ passing
through a new vertex $w$.  If $H$ is obtained from $G$ by subdividing one edge
of $G$, then $H$ is an {\it elementary subdivision} of $G$.  Edge subdivision
can reduce the edge-bandwidth considerably, but it increases the edge-bandwidth
by at most one.

\TH{11.}
If $H$ is an elementary subdivision of $G$, then
$\CL{(2B'(G)+\dlt)/3}\le B'(H)\le B'(G)+1$, where $\dlt$ is 1 if $B'(H)$ is odd
and 0 if $B'(H)$ is even, and these bounds are sharp.
\PF
Suppose that $H$ is obtained from $G$ by subdividing edge $e$.
From an optimal edge-numbering $g$ of $G$, we obtain an edge-numbering
of $H$ by augmenting the labels greater than $g(e)$ and letting the
labels of the two new edges be $g(e)$ and $g(e)+1$.  This stretches
the difference between incident labels by at most 1.

To show that this bound is sharp, compare $G=\Theta(1,2,\ldots,2)$
and $G'=\Theta(1,3,\ldots,3)$, where each has $m$ paths with common endpoints.
In Example A, we proved that $B'(G')=\CL{3(m-1)/2}$.  In $G$, let the $i$th
path have edges $a_i,b_i$ for $i<m$, with $e$ the extra edge.  The ordering
$\VEC a1{m-1},e,\VEC b1{m-1}$ yields $B'(G)\le m$.  The graph $G'$ is obtained
from $G$ by a sequence of $m-1$ elementary subdivisions, roughly half of which
must increase the edge-bandwidth.  The desired graph $H$ is the first where
the bandwidth is $m+1$.

To prove the lower bound on $B'(H)$, we consider an optimal edge-numbering
$f$ of $H$ and obtain an edge-numbering of $G$.  For the edges $e',e''$
introduced to form $H$ after deleting $e$, let $p=f(e')$ and $q=f(e'')$.
We may assume that $p<q$.  Let $r=\FL{(p+q)/2}$.  Define $g$ by leaving the
$f$-labels below $p$ and in $[r+1,q-1]$ unchanged, decreasing those in
$[p+1,r]$ and above $q$ by one, and setting $g(e)=r$.  The differences between
labels on edges belonging to both $G$ and $H$ change by at most one and increase
only when the difference is less than $B'(f)$.  For incidences involving $e$,
the incident edge $\eps$ was incident in $H$ to $e'$ or $e''$.  The difference
$|g(e)-g(\eps)|$ exceeds $B'(f)$ only if $g(\eps)<p$ or $g(\eps)>q$.  In the
first case, the difference increases by $r-p=\FL{(q-p)/2}$.  In the second,
it increases by $q-r-1=\CL{(q-p)/2}-1$.  We obtain
$B'(G)\le B'(H)+\FFR{q-p}2\le \FFR{3B'(H)}2$.  Whether $B'(H)$ is even or
odd, this establishes the bound claimed.

To show that this bound is sharp, compare $G=\Theta(1,3\ldots,3)$ and
$H=\Theta(2,3\ldots,3)$.  In $H$ let the $i$th path have edges $a_i,b_i,c_i$
for $i<m$, with $d,e$ the remaining path.  The ordering
$\VEC a1{m-1},d,\VEC b1{m-1},e,\VEC c1{m-1}$ yields $B'(H)\le m$.  
From Example A, $B'(G)=\CL{3(m-1)/2}$.  Whether $m$ is odd or even,
this example achieves the lower bound on $B'(H)$.
\qed


{\it Contracting} an edge $uv$ means deleting the edge and replacing
its endpoints by a single combined vertex $w$ inheriting all other edge
incidences involving $u$ and $v$.  Contraction tends to make a graph
denser and thus increase edge-bandwidth.  In some applications, one
restricts attention to simple graphs and thus discards loops or multiple
edges that arise under contraction.  Such a convention can discard many
edges and thus lead to a decrease in edge-bandwidth.  In particular,
contracting an edge of a clique would yield a smaller clique under this
model and thus smaller edge-bandwidth.

\smallskip
For the next result, we say that $H$ is an {\it elementary contraction} of $G$
if $H$ is obtained from $G$ by contracting one edge and keeping all other edges,
regardless of whether loops or multiple edges arise.  Edge-bandwidth is a valid
parameter for multigraphs.
\eject

\TH{12.}
If $H$ is an elementary contraction of $G$, then
$B'(G)-1\le B'(H)\le 2B'(G)-1$, and these bounds are sharp for each
value of $B'(G)$.
\PF
Let $e$ be the edge contracted to produce $H$.  For the upper bound, let $g$ be
an optimal edge-numbering of $G$, and let $f$ be the edge-numbering of $H$
produced by deleting $e$ from the numbering.  In particular, leave the
$g$-labels below $g(e)$ unchanged and decrement those above $g(e)$ by 1.
Edges incident in $H$ have distance at most two in $L(G)$, and their
distance in $L(G)$ is two only if $e$ lies between them.  Thus the
difference between their $g$-labels is at most $2B'(g)$, with equality
only if the difference between their $f$-labels is $2B'(G)-1$.

Equality holds when $G$ is the double-star (the caterpillar with
two vertices of degree $k+1$ and $2k$ vertices of degree 1) and $e$
is the central edge of $G$, so $H$ is the star $K_{1,2k}$.  We have observed
that $B'(G)=k$ and $B'(H)=2k-1$.

For the lower bound, let $f$ be an optimal edge-numbering of $H$, and let
$g$ be the edge-numbering of $G$ produced by inserting $e$ into the numbering
just above the edge $e'$ with lowest $f$-label among those incident to the
contracted vertex $w$ in $H$.  In particular, leave $f$-labels up to
$f(e')$ unchanged, augment those above $f(e')$ by 1, and let $g(e)=f(e')+1$.
The construction and the argument depend on the preservation of loops and
multiple edges.  Edges other than $e$ that are incident in $G$ are also incident
in $H$, and the difference between their labels under $g$ is at most one more
than the difference under $f$.  Edges incident to $e$ in $G$ are incident to
$e'$ in $H$ and thus have $f$-label at most $f(e')+B'(f)$.  Thus their $g$-label
differs from that of $e'$ by at most $B'(f)$.

The lower bound must be sharp for each value of $B'(G)$, because successive
contractions eventually eliminate all edges and thus reduce the bandwidth.
\qed

\SH
{5. EDGE-BANDWIDTH OF CLIQUES AND BICLIQUES}
We have computed edge-bandwidth for caterpillars and other sparse graphs.
In this section we compute edge-bandwidth for classical dense families, the
cliques and equipartite complete bipartite graphs.  Give the difficulty
of bandwidth computations, the existence of exact formulas is of as much
interest as the formulas themselves.

\TH{13.}
$B'(K_n) = \FL{n^2/4}+\CL{n/2}-2$.
\PF
{\it Lower bound.}  Consider an optimal numbering.
Among the lowest $\CH{\CL{n/2}-1}2+1$ values there must be edges
involving at least $\CL{n/2}$ vertices of $K_n$.
Among the highest $\CH{\FL{n/2}}2+1$ values there must be edges
involving at least $\FL{n/2}+1$ vertices of $K_n$.
Since $\CL{n/2}+\FL{n/2}+1>n$, some vertex has incident edges with
labels among the lowest $\CH{\CL{n/2}-1}2+1$ and among the highest
$\CH{\FL{n/2}}2+1$.  Therefore,
$$\eqalign{B'(K_n)&\ge \left[\CH n2-\CH{\FL{n/2}}2\right]-
                       \left[\CH{\CL{n/2}-1}2+1\right]\cr
       &=(\CFR n2-1)(\FFR n2)+n-1-1\cr
       &=\FFR{n^2}4+\CFR n2-2\cr}$$

{\it Upper bound.}  To achieve the bound above, let $X,Y$ be the vertex
partition with $X=\{1,\ldots,\CL{n/2}\}$ and $Y=\{\CL{n/2}+1,\ldots,n\}$.
We assign the lowest $\CH{\CL{n/2}}2$ values to the edges within $X$.
We use reverse lexicographic order, listing first the edges with
higher vertex 2, then higher vertex 3, etc.  We assign the highest
$\CH{\FL{n/2}}2$ values to the edges within $Y$ by the symmetric procedure.
Thus
$$\matrix{
      u&1&1&2&1&2&3&\cdots&\cdots&n-3&n-3&n-3&n-2&n-2&n-1&\cr
      v&2&3&3&4&4&4&\cdots&\cdots&n-2&n-1& n &n-1& n & n &\cr
  f(uv)&1&2&3&4&5&6&\cdots&\cdots&\CH n2-5& & & & &\CH n2&\cr}
$$
Note that the lowest label on an edge incident to vertex $\CL{n/2}$ is
$1+\CH{\CL{n/2}-1}2$.

The labels between these ranges are assigned to the ``cross-edges'' between $X$
and $Y$.  The cross-edges involving the vertex $\CL{n/2}\in X$ receive the
highest of the central labels, and the cross-edges involving $\CL{n/2}+1\in Y$
(but not $\CL{n/2}$) receive the lowest of these labels.  Since the highest
cross-edge label is $\CH n2-\CH{\FL{n/2}}2$ and the lowest label of an edge
incident to $\CL{n/2}$ is $1+\CH{\CL{n/2}-1}2$, the maximum difference between
labels on edges incident to $\CL{n/2}$ is precisely the lower bound on $B'(K_n)$
computed above.  This observation holds symmetrically for the edges incident to
$\CL{n/2}+1$.
\def\on#1#2{{#1\atop#2}}
$$
\on12
\on13 \on23
\on14 \on24 \on34 \
\on51 \on52 \on53 
\on61 \on62
\on71
\on81
\on72 \on82
\on63 \on73 \on83 
\on54 \on64 \on74 \on84 \
\on56 \on57 \on58
\on67 \on68
\on78
$$

We now procede iteratively.  On the high end of the remaining gap, we assign the
values to the remaining edges incident to $\CL{n/2}-1$.  Then on the low
end, we assign values to the remaining edges incident to $\CL{n/2}+2$.  We
continue alternating between the top and the bottom, completing the edges
incident to the more extreme labels as we approach the center of the numbering.
We have illustrated the resulting order for $K_8$.  Each time we insert the
remaining edges incident to a vertex of $X$, the rightmost extreme moves toward
the center at least as much from the previous extreme as the leftmost extreme
moves toward the left.  Thus the bound on the difference is maintained for the
edges incident to each vertex.  The observation is symmetric for edges incident
to vertices of $Y$.  \qed

For equipartite complete bipartite graphs, we have a similar construction
involving low vertices, high vertices, and cross-edges.

\TH{14.}
$B'(K_{n,n})=\CH{n+1}2-1$.
\PF
{\it Lower bound.}
We use the boundary bound of Proposition 3 with $k=\FL{n^2/4}+1$.
Every set of $k$ edges is together incident to at least $n+1$ vertices,
since a bipartite graph with $n$ vertices has at most $k-1$ edges.
Since $K_{n,n}$ has $2n$ vertices, at most $\FL{(n-1)^2/4}$ edges
remain when these vertices are deleted.  Thus when $\C F=k$, we have
$$
B'(K_{n,n})\ge\C{\partial(F)}\ge n^2-\FFR{(n-1)^2}4-\FFR{n^2}4-1=\CH{n+1}2-1.
$$

We construct an ordering achieving this bound.  Let $X=\{\VEC x1n\}$
and $Y=\{\VEC y1n\}$ be the partite sets.  Order the vertices as
$L=x_1,y_1,\ldots,x_n,y_n$.  We alternately finish a vertex from the beginning
of $L$ and a vertex from the end.  When finishing a vertex from the beginning,
we place its incident edges to vertices earlier in $L$ at the end of the initial
portion of the numbering $f$ that has already been determined.  When finishing a
vertex from the end of $L$, we place its incident edges to vertices later in $L$
at the beginning of the terminal portion of $f$ that has been determined.  We do
not place an edge twice.  When we have finished each vertex in each direction,
we have placed all edges in the numbering.  For example, this produces the
following edge ordering for $K_{6,6}$:
$$
\on{X}{Y}\quad 
\on11 \on21 \on12 \on22 \ 
\on31 \on32 \on13 \on23 \on33 \ 
\on41 \on42 \on43 \on14 \on24 
\on51 \on15 \on16 \on61 \ 
\on25 \on26 \on52 \on62
\on34 \on35 \on36 \on53 \on63 \ 
\on44 \on45 \on46 \on54 \on64 \ 
\on55 \on56 \on65 \on66
$$

It suffices to show that for the $j$th vertex $v_j\in L$, there are at least
$n^2-\CH{n+1}2=\CH n2$ edges that come before the first edge incident to
$v$ or after the last edge incident to $v$.  For $j=n+1$, there are exactly
$\FL{n^2/4}$ edges before the first appearance of $v_j$ and exactly
$\FL{(n-1)^2/4}$ edges after its last appearance, which matches the argument
in the lower bound.  As $j$ decreases, the leftmost appearance of $v_j$
moves leftward no more quickly than the rightmost appearance; we omit
the numerical details.  The symmetric argument applies for $j\ge n$.  \qed

\SH
{References}
\frenchspacing

\BP [1]
S.F. Assmann, G.W. Peck, M.M. Sys\l o, and J. Zak, The bandwidth of
caterpillars with hairs of length $1$ and $2$, \SIAD\ 2(1981), 387--393.
\BP [2]
P.Z.~Chinn, J.~Chv\'atalov\'a, A.K.~Dewdney, and N.E.~Gibbs, The bandwidth
problem for graphs and matrices - a survey, \JGT\ 6(1982), 223--254.
\BP [3]
F.R.K.~Chung, Labelings of graphs, in {\it Selected Topics in Graph Theory, III}
(L.~Beineke and R.~Wilson, eds.), (Academic Press 1988), 151--168.
\BP [4]
J. Chv\'atalov\'a, Optimal labelling of a product of two paths, \DM\ 11 (1975),
249--253.
\BP [5]
J. Chv\'atalov\'a and  J.Opatrn\'y, The bandwidth problem and operations on
graphs, \DM\ 61 (1986), 141--150. 
\BP [6]
J. Chv\'atalov\'a and J. Opatrn\'y, The bandwidth of theta graphs, 
{\it Utilitas Math.} 33 (1988), 9--22.
\BP [7]
D. Eichhorn, D. Mubayi, K. O'Bryant, and D.B. West,
The edge-bandwidth of theta graphs (in preparation).
\BP [8]
M.R.~Garey, R.L.~Graham, D.S.~Johnson, and D.E.~Knuth, Complexity results for
bandwidth minimization.  \SIAP\ 34(1978), 477--495.
\BP [9]
L.H. Harper, Optimal assignments of numbers to vertices, {\it J. Soc. Indust.
Appl. Math.} 12(1964), 131--135.
\BP [10]
R. Hochberg, C. McDiarmid, and M. Saks, On the bandwidth of triangulated
triangles, (Proc. 14th Brit. Comb. Conf. - Keele, 1993), \DM\ 138 (1995),
261--265.
\BP [11]
L.T.Q. Hung, M.M. Sys\l o, M.L. Weaver, and D.B. West,
Bandwidth and density for block graphs, \DM\ 189 (1998), 163--176.
\BP [12]
D.J.~Kleitman and R.V.~Vohra, Computing the bandwidth of interval graphs,
\SIDM\ 3(1990), 373--375.
\BP [13]
J. H. Mai, The bandwidth of the graph formed by $n$ meridian lines on a sphere
(Chinese, English summary), {\it J. Math. Res. Exposition} 3 (1983), 55--60.
\BP [14]
Z.~Miller, The bandwidth of caterpillar graphs, {\it Proc.\ Southeastern
Conf.{}}, {\it Congressus Numerantium} 33(1981), 235--252.
\BP [15]
H.S.~Moghadam, {\it Compression operators and a solution to the bandwidth
problem of the product of $n$ paths}, PhD Thesis, University of
California---Riverside (1983).
\BP [16]
B.~Monien, The bandwidth minimization problem for caterpillars with hair
length 3 is NP-complete, \SIAD\ {\bf7}(1986), 505--512.
\BP [17]
C.H.~Papadimitriou, the NP-completeness of the bandwidth minimization problem.
{\it Computing} 16(1976), 263--270.
\BP [18]
G.W. Peck and A. Shastri, Bandwidth of theta graphs with short paths.
\DM\ 103 (1992), 177--187. 
\BP [19]
L. Smithline, Bandwidth of the complete $k$-ary tree,
\DM\ 142(1995), 203--212.
\BP [20]
A.P.~Sprague, An $O(n\log n)$ algorithm for bandwidth of interval graphs,
\SIDM\  {\bf7}(1994), 213--220.
\BP [21]
M.M.~Sys\l o and J.~Zak, The bandwidth problem: critical subgraphs and the
solution for caterpillars, \ADM\ 16(1982), 281--286.  See also Comp.\ Sci.\
Dept.\ Report CS-80-065, Washington State Univ.\ (1980).
\BP [22]
J.F. Wang, D.B. West, and B. Yao, Maximum bandwidth under edge addition,
\JGT\ 20(1995), 87--90. 
\bye